# Explicit inversion of cosh-weighted Finite Hilbert Transform

**Jiangsheng You**

Cubic Imaging LLC, Andover, MA 01810

jshyou@gmail.com

## Abstract

A contour integral approach is presented to rederive the Bertola-Katsevich-Tovbis inversion formulas of the cosh-weighted finite Hilbert Transform without reducing to a Riemann-Hilbert problem. The contour integral approach yields the explicit expression of the cosh-weighted Hilbert transforms of two classes of exponential Chebyshev functions. Four numerical experiments are performed to show the superior numerical stability of the Bertola-Katsevich-Tovbis inversion formulas for large exponential constants from using the sine and cosine transforms at Chebyshev nodes.

## I. Mathematical preliminary

The finite Hilbert transform (FHT) over the open unit interval $I=(-1,\ 1)$ is defined by

$$F(s)=Hf(s)=\frac{1}{\pi}\int_{-1}^{1}\frac{f(t)}{s-t}dt\ . \tag{1.1}$$

For a complex constant $\mu\in C$, the cosh-weighted FHT (chFHT) is defined as

$$F_{\mu}(s)=H_{\mu}f(s)=\frac{1}{\pi}\int_{-1}^{1}\frac{\cosh[\mu(s-t)]}{s-t}f(t)dt\ . \tag{1.2}$$

Throughout this paper we consider $F(s)$ and $F_{\mu}(s)$ in $I$ though they can be defined on $R^1$. The inversion of FHT has been well established in literature [1-4, 8, 12, 14]. The applications of (1.1) are referenced to [6-7, 10, 21]. The study of (1.2) was initiated by [13] in which the reconstruction of half scan data from exponential Radon transform was converted into the inversion of (1.2). The early numerical inversions of (1.2) were studied in [5, 11, 14] but the analysis of numerical stability was very limited. The author studied the properties of (1.2) in weighted $L^2$ space [18-19] and continued the research on the inversion of (1.2) following the breakthrough idea of [2] in which the range condition and explicit inversions of (1.2) were derived by reducing to a vector Riemann-Hilbert problem (RHP). In this paper, the author presents a direct contour integral approach to rederive the Bertola-Katsevich-Tovbis inversion formulas by only using the Sokhotski-Plemelj and Poincaré-Bertrand formulas. For two cases of exponential Chebyshev functions studied in [13, 20], the explicit chFHTs are obtained by the contour integral method. Four numerical experiments are performed to show the superior numerical stability of the Bertola-Katsevich-Tovbis inversion formulas for a large value range of $\mu$ and in the presence of noise by using the sine and cosine transforms at Chebyshev nodes.

For $1\le p<\infty$, $L^p(I)$ is the $p$-integrable Banach space. The FHT $F(s)$ is continuous in $L^p(I)$, $1<p<\infty$ [15]. The author studied the chFHT in the weighted $L^2(I)$ in [18-19] taking the advantages of numerical analysis at Chebyshev nodes. We define two weighted spaces

$$L_m^2(I)=\{f(t):\int_{-1}^{1}|f(t)|^2\sqrt{1-t^2}dt<\infty\}\ , \tag{1.3}$$

$$L_d^2(I)=\{f(t):\int_{-1}^{1}|f(t)|^2\frac{1}{\sqrt{1-t^2}}dt<\infty\}. \tag{1.4}$$

Subscriptions $m$ and $d$ stand for the multiplication and division of $\sqrt{1-t^2}$, respectively. Notice that $L_d^2(I)\subset L^2(I)\subset L_m^2(I)$ and $C_0^\infty(I)$ is dense in these linear spaces. The weighted spaces $L_m^2(I)$ and $L_d^2(I)$ are Hilbert spaces with the following inner products

$$<f,g>_m=\frac{1}{\pi}\int_{-1}^{1}f(t)g^*(t)\sqrt{1-t^2}\,dt,\qquad <f,g>_d=\frac{1}{\pi}\int_{-1}^{1}\frac{f(t)g^*(t)}{\sqrt{1-t^2}}dt. \tag{1.5}$$

Here the star stands for the complex conjugate. Notice that $\|1/\sqrt{1-t^2}\|_{L_m^2}=1$ but $1/\sqrt{1-t^2}\notin L^2(I)$. In this paper, motivated by [2] and earlier studies [18-19], the author presents a unified contour integral approach to rederive the range and inversion of (1.2) in $L_m^2(I)$ and $L_d^2(I)$. The proposed approach only uses the Poincaré-Bertrand and Sokhotski-Plemelj formulas.

For the Poincaré-Bertrand formula, we adopt the expression given in [8] to handle the term $\cosh[\mu(s-t)]$ more conveniently. For a function $\varphi(u,s)$, the formula reads

$$\frac{1}{\pi}\int_{-1}^{1}\frac{1}{s-t}[\frac{1}{\pi}\int_{-1}^{1}\frac{\varphi(u,s)}{s-u}du]ds=\varphi(t,t)+\frac{1}{\pi}\int_{-1}^{1}\frac{du}{t-u}\frac{1}{\pi}\int_{-1}^{1}[\frac{1}{t-s}-\frac{1}{u-s}]\varphi(u,s)ds. \tag{1.6}$$

For the Sokhotski-Plemelj formula, we use the following expression to consider the jump cross the boundary. The formula reads

$$\lim_{\substack{\varepsilon>0\\ \varepsilon\to 0}}\frac{1}{\pi}\int_{-1}^{1}\frac{\phi(s\pm i\varepsilon)}{t-(s\pm i\varepsilon)}ds=\pm i\phi^{\pm}(t)+\frac{1}{\pi}\int_{-1}^{1}\frac{1}{t-s}\phi^{\pm}(s)ds. \tag{1.7}$$

Here $\phi^{\pm}(s)=\lim_{\varepsilon>0,\varepsilon\to 0}\phi(s\pm i\varepsilon)$, and the detailed exposition of the Sokhotski-Plemelj formula can be found in [8].

## II. Identities and inversions of chFHT

The arguments of [2] are to cast the inversion of (1.2) to a vector RHP to derive the left and right inverses, and then use two lemmas to compute the inversion formula. Several main results such as (3.10), (4.6) and Theorem B.4 of [2] are derived through the vector RHP while Theorem 6.2 of [2] was obtained by lengthy computations. This paper uses (1.6) and (1.7) to rederive these Bertola-Katsevich-Tovbis formulas through the contour integrals of a test function in the infinity and around $[-1,\ 1]$. The structure of this section is to derive five identities in Lemma 1 and then use them to obtain two main formulas in Theorem 1 and a family of formulas in Remark 3.

**Lemma 1**. *For* $t,u\in(-1,\ 1)$ *and* $\mu\in C$*, we have the following identities*

$$\frac{1}{\pi}\int_{-1}^{1}\cosh[\mu(u-s)]\frac{\cos(\mu\sqrt{1-s^2})}{\sqrt{1-s^2}}ds=\cosh(\mu u), \tag{2.1}$$

$$\frac{1}{\pi}\int_{-1}^{1}\cosh[\mu(u-s)]\frac{\cos(\mu\sqrt{1-s^2})}{\sqrt{1-s^2}}s\,ds=-\frac{\mu}{2}\sinh(\mu u), \tag{2.2}$$

$$\frac{1}{\pi}\int_{-1}^{1}\frac{\cosh[\mu(u-s)]}{t-s}\sin(\mu\sqrt{1-s^2})ds=\sinh(\mu u)-\sinh[\mu(u-t)]\cos(\mu\sqrt{1-t^2}), \tag{2.3}$$

$$\frac{1}{\pi}\int_{-1}^{1}\frac{\cosh[\mu(u-s)]}{t-s}\frac{\cos(\mu\sqrt{1-s^2})}{\sqrt{1-s^2}}ds=\sinh[\mu(u-t)]\sin(\mu\sqrt{1-t^2})\frac{1}{\sqrt{1-t^2}}, \tag{2.4}$$

$$\frac{1}{\pi}\int_{-1}^{1}\frac{\cosh[\mu(u-s)]}{t-s}\cos(\mu\sqrt{1-s^2})\sqrt{1-s^2}\,ds \\ = \sinh[\mu(u-t)]\sin(\mu\sqrt{1-t^2})\sqrt{1-t^2}+t\cosh(\mu u)-\frac{\mu}{2}\sinh(\mu u). \tag{2.5}$$

**Proof**. Consider a complex test function

$$\Phi_\mu(u,z)=\frac{1}{2}e^{\mu(u-z+\sqrt{z^2-1})},\ u\in(-1,\ 1). \tag{2.6}$$

Choose the cuts of $\sqrt{z^2-1}$ satisfying

$$\lim_{|z|\to\infty}\frac{\sqrt{z^2-1}}{z}=1,\ \lim_{\varepsilon>0,\varepsilon\to 0}\sqrt{(s\pm i\varepsilon)^2-1}=\pm i\sqrt{1-s^2}\ \text{ for } s\in(-1,\ 1). \tag{2.7}$$

With such selection of cuts, $\Phi_\mu(u,z)$ is analytical in $C\setminus[-1,\ 1]$ and

$$\Phi_\mu(u,\infty)=\lim_{|z|\to\infty}\Phi_\mu(u,z)=\frac{1}{2}e^{\mu u},\quad \lim_{|z|\to\infty}[\Phi_\mu(u,z)-\Phi_\mu(u,\infty)]z=-\frac{1}{4}\mu e^{\mu u}, \tag{2.8}$$

$$\Phi_\mu^{\pm}(u,s)=\lim_{\varepsilon>0,\varepsilon\to 0}\Phi_\mu(u,s\pm i\varepsilon)=\frac{1}{2}e^{\mu(u-s\pm i\sqrt{1-s^2})},\ s\in(-1,\ 1). \tag{2.9}$$

To be consistent with the singular integrals around $[-1,\ 1]$, the contour integrals are clockwise on large circles so that $\oint_{|z|=d}(1/z)dz=-2\pi i$. Here we only present the computations for (2.3) and (2.4). For $|z|\to\infty$, it is straightforward to derive

$$\frac{1}{2\pi i}\int_{|z|=d,d>1}\frac{1}{t-z}[\Phi_\mu(u,z)-\Phi_{-\mu}(u,z)]dz=\lim_{d\to\infty}\frac{1}{2\pi i}\oint_{|z|=d}\frac{\Phi_\mu(u,z)-\Phi_{-\mu}(u,z)}{t-z}dz=\sinh(\mu u), \tag{2.10}$$

$$\frac{1}{2\pi i}\int_{|z|=d,d>1}\frac{1}{t-z}\frac{\Phi_\mu(u,z)+\Phi_{-\mu}(u,z)}{\sqrt{z^2-1}}dz=\lim_{d\to\infty}\frac{1}{2\pi i}\oint_{|z|=d}\frac{1}{t-z}\frac{\Phi_\mu(u,z)+\Phi_{-\mu}(u,z)}{\sqrt{z^2-1}}dz=0. \tag{2.11}$$

Around $[-1,\ 1]$, the left side of (2.10) becomes

$$\begin{aligned}
&\frac{1}{2\pi i}\int_{|z|=d,d>1}\frac{1}{t-z}[\Phi_\mu(u,z)-\Phi_{-\mu}(u,z)]dz\\
&=\lim_{\substack{\varepsilon>0\\ \varepsilon\to 0}}\frac{1}{2\pi i}\int_{-1}^{1}\frac{1}{t-(s+i\varepsilon)}[\Phi_\mu(u,s+i\varepsilon)-\Phi_{-\mu}(u,s+i\varepsilon)]ds\\
&-\lim_{\substack{\varepsilon>0\\ \varepsilon\to 0}}\frac{1}{2\pi i}\int_{-1}^{1}\frac{1}{t-(s-i\varepsilon)}[\Phi_\mu(u,s-i\varepsilon)-\Phi_{-\mu}(u,s-i\varepsilon)]ds\\
&=\sinh[\mu(u-t)]\cosh(i\mu\sqrt{1-t^2})+\frac{1}{\pi i}\int_{-1}^{1}\frac{\cosh[\mu(u-s)]}{t-s}\sinh(i\mu\sqrt{1-s^2})ds\\
&=\sinh[\mu(u-t)]\cos(\mu\sqrt{1-t^2})+\frac{1}{\pi}\int_{-1}^{1}\frac{\cosh[\mu(u-s)]}{t-s}\sin(\mu\sqrt{1-s^2})ds.
\end{aligned} \tag{2.12}$$

Here the first limit of (2.12) is

$$\frac{1}{2}[\Phi_\mu^{+}(u,t)-\Phi_{-\mu}^{+}(u,t)]+\frac{1}{2\pi i}\int_{-1}^{1}\frac{1}{t-s}[\Phi_\mu^{+}(u,s)-\Phi_{-\mu}^{+}(u,s)]ds; \tag{2.13}$$

and the second limit of (2.12) is

$$-\frac{1}{2}[\Phi_\mu^{-}(u,t)-\Phi_{-\mu}^{-}(u,t)]+\frac{1}{2\pi i}\int_{-1}^{1}\frac{1}{t-s}[\Phi_\mu^{-}(u,s)-\Phi_{-\mu}^{-}(u,s)]ds. \tag{2.14}$$

Then (2.3) is derived by (2.10) and (2.12). Around $[-1,\ 1]$, the left side of (2.11) becomes

$$
\begin{aligned}
&\frac{1}{2\pi i}\int_{|z|=d,d>1}\frac{1}{t-z}\frac{\Phi_{\mu}(u,z)+\Phi_{-\mu}(u,z)}{\sqrt{z^2-1}}dz\\
&=\lim_{\substack{\varepsilon>0\\ \varepsilon\to 0}}\frac{1}{2\pi i}\int_{-1}^{1}\frac{1}{t-(s+i\varepsilon)}\frac{\Phi_{\mu}(u,s+i\varepsilon)+\Phi_{-\mu}(u,s+i\varepsilon)}{\sqrt{(s+i\varepsilon)^2-1}}ds\\
&-\lim_{\substack{\varepsilon>0\\ \varepsilon\to 0}}\frac{1}{2\pi i}\int_{-1}^{1}\frac{1}{t-(s-i\varepsilon)}\frac{\Phi_{\mu}(u,s-i\varepsilon)+\Phi_{-\mu}(u,s-i\varepsilon)}{\sqrt{(s-i\varepsilon)^2-1}}ds \qquad (2.15)\\
&=\frac{\sinh[\mu(u-t)]\sinh(i\mu\sqrt{1-t^2})}{i\sqrt{1-t^2}}+\frac{1}{\pi i}\int_{-1}^{1}\frac{\cosh[\mu(s-u)]}{t-s}\frac{\cosh(i\mu\sqrt{1-s^2})}{i\sqrt{1-s^2}}ds\\
&=\frac{\sinh[\mu(u-t)]\sin(\mu\sqrt{1-t^2})}{\sqrt{1-t^2}}-\frac{1}{\pi}\int_{-1}^{1}\frac{\cosh[\mu(s-u)]}{t-s}\frac{\cos(\mu\sqrt{1-s^2})}{\sqrt{1-s^2}}ds.
\end{aligned}
$$

Here the first limit of (2.15) is

$$
\frac{1}{2}\frac{\Phi_{\mu}^{+}(u,t)+\Phi_{-\mu}^{+}(u,t)}{i\sqrt{1-t^2}}+\frac{1}{2\pi i}\int_{-1}^{1}\frac{1}{t-s}\frac{\Phi_{\mu}^{+}(u,s)+\Phi_{-\mu}^{+}(u,s)}{i\sqrt{1-s^2}}ds, \qquad (2.16)
$$

and the second limit of (2.15) is

$$
\frac{1}{2}\frac{\Phi_{\mu}^{-}(u,t)+\Phi_{-\mu}^{-}(u,t)}{i\sqrt{1-t^2}}+\frac{1}{2\pi i}\int_{-1}^{1}\frac{1}{t-s}\frac{\Phi_{\mu}^{-}(u,s)+\Phi_{-\mu}^{-}(u,s)}{-i\sqrt{1-s^2}}ds. \qquad (2.17)
$$

It follows that (2.4) is derived by (2.11) and (2.15). The derivations of (2.1) and (2.2) follow the same argument. Identity (2.5) can be derived by using (2.1), (2.2) and (2.4) as follows

$$
\begin{aligned}
&\frac{1}{\pi}\int_{-1}^{1}\frac{\cosh[\mu(u-s)]}{t-s}\cos(\mu\sqrt{1-s^2})\sqrt{1-s^2}ds\\
&=\frac{1}{\pi}\int_{-1}^{1}\frac{\cosh[\mu(u-s)]}{t-s}\frac{\cos(\mu\sqrt{1-s^2})}{\sqrt{1-s^2}}(1-s^2)ds\\
&=\frac{1}{\pi}\int_{-1}^{1}\frac{\cosh[\mu(u-s)]}{t-s}\frac{\cos(\mu\sqrt{1-s^2})}{\sqrt{1-s^2}}(1-t^2+t^2-s^2)ds \qquad (2.18)\\
&=\sinh[\mu(u-t)]\sin(\mu\sqrt{1-t^2})\sqrt{1-t^2}+\frac{1}{\pi}\int_{-1}^{1}\cosh[\mu(u-s)]\frac{\cos(\mu\sqrt{1-s^2})}{\sqrt{1-s^2}}(t+s)ds\\
&=\sinh[\mu(u-t)]\sin(\mu\sqrt{1-t^2})\sqrt{1-t^2}+t\cosh(\mu u)-\frac{\mu}{2}\sinh(\mu u).
\end{aligned}
$$

The proof of (2.1-2.5) is complete. □

**Remark 1**. Identity (2.3) is equivalent to Lemma A.2 of [2] and (2.5) is equivalent to Lemma A.1 of [2]. Taking $u=t$, (2.4) yields Theorem 3.3 of [2] in which $\cos(\mu\sqrt{1-t^2})/\sqrt{1-t^2}$ is the non-trivial function in the null space of $H_\mu$ in $L_m^2(I)$.

The test function $\Phi_\mu(s,z)$ of (2.6) is very crucial in all the computations of Lemma A.1 and Lemma A.2 in [2] and is also the key factor to derive all the identities in Lemma 1. Next, we use Lemma 1 to prove the Bertola-Katsevich-Tovbis inversion formulas.

**Theorem 1**. *For $\mu\in C$, we rephrase the Bertola-Katsevich-Tovbis inversion formulas of chFHT in two identities. For $f(t)\in L_d^2(I)$, we have*

$$f(t)=\cos(\mu\sqrt{1-t^2})\sqrt{1-t^2}[\frac{1}{\pi}\int_{-1}^{1}\frac{F_\mu(s)}{s-t}\frac{\cos(\mu\sqrt{1-s^2})}{\sqrt{1-s^2}}ds] + \sin(\mu\sqrt{1-t^2})[\frac{1}{\pi}\int_{-1}^{1}\frac{F_\mu(s)}{s-t}\sin(\mu\sqrt{1-s^2})ds]. \tag{2.19}$$

*For* $f(t)\in L_m^2(I)$*, we have*

$$f(t)=\frac{\cos(\mu\sqrt{1-t^2})}{\sqrt{1-t^2}}[\frac{1}{\pi}\int_{-1}^{1}\frac{F_\mu(s)}{s-t}\cos(\mu\sqrt{1-s^2})\sqrt{1-s^2}ds+\frac{1}{\pi}\int_{-1}^{1}\cosh(\mu t)f(t)dt] + \sin(\mu\sqrt{1-t^2})[\frac{1}{\pi}\int_{-1}^{1}\frac{F_\mu(s)}{s-t}\sin(\mu\sqrt{1-s^2})ds]. \tag{2.20}$$

**Proof**. We compute $\sin(\mu\sqrt{1-t^2})[\frac{1}{\pi}\int_{-1}^{1}\frac{F_\mu(s)}{s-t}\sin(\mu\sqrt{1-s^2})ds]$ as follows

$$\begin{aligned}&\sin(\mu\sqrt{1-t^2})\frac{1}{\pi}\int_{-1}^{1}\frac{F_\mu(s)}{s-t}\sin(\mu\sqrt{1-s^2})ds\\ &=\sin(\mu\sqrt{1-t^2})\frac{1}{\pi}\int_{-1}^{1}\frac{\sin(\mu\sqrt{1-s^2})}{s-t}[\frac{1}{\pi}\int_{-1}^{1}\frac{\cosh[\mu(s-u)]}{s-u}f(u)du]ds\\ &=\sin(\mu\sqrt{1-t^2})\frac{1}{\pi}\int_{-1}^{1}\frac{1}{s-t}[\frac{1}{\pi}\int_{-1}^{1}\frac{\cosh[\mu(s-u)]\sin(\mu\sqrt{1-s^2})f(u)}{s-u}du]ds.\end{aligned} \tag{2.21}$$

Apply (1.6) to (2.21) with $\varphi(u,s)=\cosh[\mu(u-s)\sin(\mu\sqrt{1-s^2})f(u)$, (2.3) yields

$$\begin{aligned}&\sin(\mu\sqrt{1-t^2})\frac{1}{\pi}\int_{-1}^{1}\frac{F_\mu(s)}{s-t}\sin(\mu\sqrt{1-s^2})ds\\ &=\sin^2(\mu\sqrt{1-t^2})f(t)+\sin(\mu\sqrt{1-t^2})\frac{1}{\pi}\int_{-1}^{1}\frac{1}{t-u}[\frac{1}{\pi}\int_{-1}^{1}(\frac{1}{t-s}-\frac{1}{u-s})\varphi(u,s)ds]du\\ &=\sin^2(\mu\sqrt{1-t^2})f(t)-\sin(\mu\sqrt{1-t^2})\cos(\mu\sqrt{1-t^2})\int_{-1}^{1}\frac{\sinh[\mu(u-t)]}{t-u}f(u)du.\end{aligned} \tag{2.22}$$

We compute $\cos(\mu\sqrt{1-t^2})\sqrt{1-t^2}\,\frac{1}{\pi}\int_{-1}^{1}\frac{F_\mu(s)}{s-t}\frac{\cos(\mu\sqrt{1-s^2})}{\sqrt{1-s^2}}ds$ as follows

$$\begin{aligned}&\cos(\mu\sqrt{1-t^2})\sqrt{1-t^2}\,\frac{1}{\pi}\int_{-1}^{1}\frac{F_\mu(s)}{s-t}\frac{\cos(\mu\sqrt{1-s^2})}{\sqrt{1-s^2}}ds\\ &=\cos(\mu\sqrt{1-t^2})\sqrt{1-t^2}\,\frac{1}{\pi}\int_{-1}^{1}\frac{1}{s-t}\frac{\cos(\mu\sqrt{1-s^2})}{\sqrt{1-s^2}}[\frac{1}{\pi}\int_{-1}^{1}\frac{\cosh[\mu(s-u)]}{s-u}f(u)du]ds\\ &=\cos(\mu\sqrt{1-t^2})\sqrt{1-t^2}\,\frac{1}{\pi}\int_{-1}^{1}\frac{1}{s-t}[\frac{1}{\pi}\int_{-1}^{1}\frac{\cosh[\mu(u-s)]}{s-u}\frac{\cos(\mu\sqrt{1-s^2})}{\sqrt{1-s^2}}f(u)du]ds.\end{aligned} \tag{2.23}$$

Apply (1.6) to (2.23) with $\varphi(u,s)=\cosh[\mu(u-s)]\frac{\cos(\mu\sqrt{1-s^2})}{\sqrt{1-s^2}}f(u)$, (2.4) yields

$$\begin{aligned}
&\cos(\mu\sqrt{1-t^2})\sqrt{1-t^2}\,\frac{1}{\pi}\int_{-1}^{1}\frac{F_\mu(s)}{s-t}\frac{\cos(\mu\sqrt{1-s^2})}{\sqrt{1-s^2}}ds \\
&=\cos^2(\mu\sqrt{1-t^2})f(t)+\cos(\mu\sqrt{1-t^2})\sqrt{1-t^2}\,\frac{1}{\pi}\int_{-1}^{1}\frac{1}{t-u}\Big[\frac{1}{\pi}\int_{-1}^{1}\Big(\frac{1}{t-s}-\frac{1}{u-s}\Big)\varphi(u,s)ds\Big]du \\
&=\cos^2(\mu\sqrt{1-t^2})f(t)+\sin(\mu\sqrt{1-t^2})\cos(\mu\sqrt{1-t^2})\int_{-1}^{1}\frac{\sinh[\mu(u-t)]}{t-u}f(u)du.
\end{aligned}\tag{2.24}$$

Recall that $\cos^2 z+\sin^2 z\equiv 1$, $z\in C$, we obtain (2.19) by adding (2.22) and (2.24) together. The similar computations yield (2.20). □

**Remark 2**. Formula (2.19) was originally proved as part of Theorem 5.1 in [2] by using the RHP for $f\in L^p(I)$, $p>2$. The author wants to emphasize that the discovery of (2.19) is remarkable and due to the vector RHP method, otherwise it is very difficult to even imagine the existence of explicit inversion formulas. Formula (2.20) was derived in [2] by using complicated computations. In this paper, Lemma 1 and Theorem 1 provide one unified way to prove both inversion formulas without using RHP. If $\mu=0$, (2.19) and (2.20) degenerate to the classical inversion of FHT. In Theorem 2, we rephrase Theorem 5.1 of [2] for $L_d^2(I)$ and $L_m^2(I)$.

**Theorem 2**. *We assume* $\mu\in C$ *and define two subspaces* $\text{Ł}_m^2(I)$ *and* $\text{Ł}_d^2(I)$ *as*

$$\text{Ł}_m^2(I)=\{f(t)\in L_m^2(I):\int_{-1}^{1}\cosh(\mu t)f(t)dt=0\}\,,\tag{2.25}$$

$$\text{Ł}_d^2(I)=\{f(t)\in L_d^2(I):\int_{-1}^{1}\frac{\cos(\mu\sqrt{1-t^2})}{\sqrt{1-t^2}}f(t)dt=0\}\,.\tag{2.26}$$

*The null space of* $H_\mu$ *in* $L_d^2(I)$ *is trivial and the null space of* $H_\mu$ *in* $L_m^2(I)$ *is defined by* $\{c\cos(\mu\sqrt{1-t^2})/\sqrt{1-t^2}\}$, $c\in C$. $H_\mu$ *is continuous from* $\text{Ł}_m^2(I)$ *to* $L_m^2(I)$ *and from* $L_d^2(I)$ *to* $\text{Ł}_d^2(I)$, *respectively.*

*For* $F_\mu(s)\in\text{Ł}_d^2(I)$, *the inverse* ${}^dH_\mu^{-1}$ *is given by*

$$\begin{aligned}
{}^dH_\mu^{-1}F_\mu(t)=&\cos(\mu\sqrt{1-t^2})\sqrt{1-t^2}\,\frac{1}{\pi}\int_{-1}^{1}\frac{1}{s-t}\frac{\cos(\mu\sqrt{1-s^2})}{\sqrt{1-s^2}}F_\mu(s)ds \\
&+\sin(\mu\sqrt{1-t^2})\frac{1}{\pi}\int_{-1}^{1}\frac{1}{s-t}\sin(\mu\sqrt{1-s^2})F_\mu(s)ds,
\end{aligned}\tag{2.27}$$

*For* $F_\mu(s)\in L_m^2(I)$, *the inverse* ${}^mH_\mu^{-1}$ *is given by*

$$\begin{aligned}
{}^mH_\mu^{-1}F_\mu(t)=&\frac{\cos(\mu\sqrt{1-t^2})}{\sqrt{1-t^2}}\frac{1}{\pi}\int_{-1}^{1}\frac{1}{s-t}\cos(\mu\sqrt{1-s^2})\sqrt{1-s^2}F_\mu(s)ds \\
&+\sin(\mu\sqrt{1-t^2})\frac{1}{\pi}\int_{-1}^{1}\frac{1}{s-t}\sin(\mu\sqrt{1-s^2})F_\mu(s)ds.
\end{aligned}\tag{2.28}$$

*Then* ${}^mH_\mu^{-1}$ *is continuous from* $L_m^2(I)$ *to* $\text{Ł}_m^2(I)$ *and* ${}^dH_\mu^{-1}$ *is continuous from* $\text{Ł}_d^2(I)$ *to* $L_d^2(I)$. *The range of* $H_\mu$ *in* $\text{Ł}_m^2(I)$ *is* $L_m^2(I)$ *and the range of* $H_\mu$ *in* $L_d^2(I)$ *is* $\text{Ł}_d^2(I)$.

**Proof**. The results described here mainly rephrase Theorem 5.1 of [2] for $L_d^2(I)$ and $L_m^2(I)$, thus we omit the details except confirming that ${}^mH_\mu^{-1}(L_m^2(I))=\text{Ł}_m^2(I)$. Let $g(t)={}^mH_\mu^{-1}F_\mu(t)$, we compute $\int_{-1}^{1}\cosh(\mu t)g(t)dt$. Take $u=0$ in (2.4), the first term of (2.28) is

$$
\begin{aligned}
&\int_{-1}^{1}\cosh(\mu t)[\frac{\cos(\mu\sqrt{1-t^2})}{\sqrt{1-t^2}}\frac{1}{\pi}\int_{-1}^{1}\frac{1}{s-t}\cos(\mu\sqrt{1-s^2})\sqrt{1-s^2}F_\mu(s)ds]dt \\
&=\int_{-1}^{1}\cos(\mu\sqrt{1-s^2})\sqrt{1-s^2}F_\mu(s)[\frac{1}{\pi}\int_{-1}^{1}\frac{\cosh(\mu t)}{s-t}[\frac{\cos(\mu\sqrt{1-t^2})}{\sqrt{1-t^2}}dt]ds \\
&=-\int_{-1}^{1}\cos(\mu\sqrt{1-s^2})\sqrt{1-s^2}F_\mu(s)[\sinh(\mu s)\sin(\mu\sqrt{1-s^2})\frac{1}{\sqrt{1-s^2}}]ds \\
&=-\int_{-1}^{1}\cos(\mu\sqrt{1-s^2})\sin(\mu\sqrt{1-s^2})\sinh(\mu s)F_\mu(s)ds.
\end{aligned}
\tag{2.29}
$$

Take $u=0$ in (2.3), the second term of (2.28) is

$$
\begin{aligned}
&\int_{-1}^{1}\cosh(\mu t)[\sin(\mu\sqrt{1-t^2})\frac{1}{\pi}\int_{-1}^{1}\frac{1}{s-t}\sin(\mu\sqrt{1-s^2})F_\mu(s)ds]dt \\
&=\int_{-1}^{1}\sin(\mu\sqrt{1-s^2})F_\mu(s)[\frac{1}{\pi}\int_{-1}^{1}\frac{1}{s-t}\cosh(\mu t)\sin(\mu\sqrt{1-t^2})dt]ds \\
&=\int_{-1}^{1}\sin(\mu\sqrt{1-s^2})F_\mu(s)[\sinh(\mu s)\cos(\mu\sqrt{1-s^2})]ds \\
&=\int_{-1}^{1}\cos(\mu\sqrt{1-s^2})\sin(\mu\sqrt{1-s^2})\sinh(\mu s)F_\mu(s)ds.
\end{aligned}
\tag{2.30}
$$

Combining (2.29) and (2.30) yields $\int_{-1}^{1}\cosh(\mu t)f(t)dt=0$. This completes the proof. □

**Remark 3**. A family of two-parameter inversions are given by (B.3) in [2] through the RHP. Here we mention that the same family of inversions can be derived through Lemma 1. Combining identities from (2.2) to (2.5), it is straightforward to obtain

$$
\begin{aligned}
&\frac{1}{\pi}\int_{-1}^{1}\frac{\cosh[\mu(u-s)]}{t-s}\frac{\alpha s^2+\beta s+\gamma}{\sqrt{1-s^2}}\cos(\mu\sqrt{1-s^2})ds \\
&=\sinh[\mu(u-t)]\sin(\mu\sqrt{1-t^2})\frac{\alpha t^2+\beta t+\gamma}{\sqrt{1-t^2}}-\alpha t\cosh(\mu u)+\frac{\mu(\alpha+\beta)}{2}\sinh(\mu u).
\end{aligned}
\tag{2.31}
$$

Assume that $\alpha t^2+\beta t+\gamma\neq 0$ in $(-1,\ 1)$, lengthy derivation yields

$$
\begin{aligned}
f(t)=&\frac{\sqrt{1-t^2}}{\alpha t^2+\beta t+\gamma}\cos(\mu\sqrt{1-t^2})[\frac{1}{\pi}\int_{-1}^{1}\frac{F_\mu(s)}{s-t}\cos(\mu\sqrt{1-s^2})\frac{\alpha s^2+\beta s+\gamma}{\sqrt{1-s^2}}ds] \\
&-\frac{\alpha}{\sqrt{1-t^2}}\cos(\mu\sqrt{1-t^2})[\frac{1}{\pi}\int_{-1}^{1}\cosh(\mu t)f(t)dt] \\
&+\sin(\mu\sqrt{1-t^2})[\frac{1}{\pi}\int_{-1}^{1}\frac{F_\mu(s)}{s-t}\sin(\mu\sqrt{1-s^2})ds].
\end{aligned}
\tag{2.32}
$$

If $\alpha\neq 0$, $\int_{-1}^{1}\cosh(\mu t)f(t)dt$ is required to be known for an exact inversion by (2.32). For $\alpha=0$, $\beta=\pm 1$ and $\gamma=1$, we have

$$
\begin{aligned}
f(t)=&\frac{1\pm t}{\sqrt{1-t^2}}\cos(\mu\sqrt{1-t^2})[\frac{1}{\pi}\int_{-1}^{1}\frac{F_\mu(s)}{s-t}\cos(\mu\sqrt{1-s^2})\frac{\sqrt{1-s^2}}{1\pm s}ds] \\
&+\sin(\mu\sqrt{1-t^2})[\frac{1}{\pi}\int_{-1}^{1}\frac{F_\mu(s)}{s-t}\sin(\mu\sqrt{1-s^2})ds].
\end{aligned}
\tag{2.33}
$$

# III. Exponential Chebyshev functions

It is well known that the Chebyshev polynomials satisfy the following relations

$$\frac{1}{\pi}\int_{-1}^{1}\frac{1}{t-s}\frac{T_n(s)}{\sqrt{1-s^2}}ds=-U_{n-1}(t), \tag{3.1}$$

$$\frac{1}{\pi}\int_{-1}^{1}\frac{1}{t-s}U_{n-1}(s)\sqrt{1-s^2}ds=T_n(t). \tag{3.2}$$

Hereafter we denote by $T_n(t)$ and $U_n(t)$ the Chebyshev polynomials of the first and second kind, respectively. In this section we introduce two classes of exponential Chebyshev functions studied in [13, 20] for the inversion of the exponential Radon transform and then derive the explicit expressions of the chFHT of these exponential Chebyshev functions. For $n\geq 0$, the exponential Chebyshev functions are defined as

$$T_{\mu,n}(s)=\frac{1}{2}[e^{-i\mu\sqrt{1-s^2}}(s+i\sqrt{1-s^2})^n+e^{i\mu\sqrt{1-s^2}}(s-i\sqrt{1-s^2})^n], \tag{3.3}$$

$$U_{\mu,n}(s)=\frac{1}{2i\sqrt{1-s^2}}[e^{-i\mu\sqrt{1-s^2}}(s+i\sqrt{1-s^2})^{n+1}-e^{i\mu\sqrt{1-s^2}}(s-i\sqrt{1-s^2})^{n+1}]. \tag{3.4}$$

Notice that $T_{0,n}(t)=T_n(t)$ and $U_{0,n}(t)=U_n(t)$. Seemingly complicated (3.3-3.4) are actually simple in the trigonometric representation

$$T_{\mu,n}(\cos\theta)=\cos(n\theta-\mu\sin\theta), \tag{3.5}$$

$$U_{\mu,n}(\cos\theta)=\frac{\sin((n+1)\theta-\mu\sin\theta)}{\sin\theta}. \tag{3.6}$$

The exponential Chebyshev functions play a very important role to explicitly express the inverse exponential Radon transform in the series of circular harmonics [13, 20] for the exterior problem. By using the same contour integrals with different test functions, we are able to express the chFHTs of $T_{\mu,n}(s)$ and $U_{\mu,n}(s)$. We state the explicit formulas in Theorem 3.

**Theorem 3**. *For the exponential Chebyshev functions (3.3) and (3.4), if* $n\geq 1$ *we have*

$$\frac{1}{\pi}\int_{-1}^{1}\frac{\cosh[\mu(t-s)]}{t-s}\frac{T_{\mu,n}(s)}{\sqrt{1-s^2}}ds=-e^{-\mu t}\sum_{k=0}^{n-1}\frac{1}{k!}\mu^k U_{n-1-k}(t), \tag{3.7}$$

$$\frac{1}{\pi}\int_{-1}^{1}\frac{\cosh[\mu(t-s)]}{t-s}U_{\mu,n-1}(s)\sqrt{1-s^2}ds=\frac{e^{-\mu t}}{2}\begin{cases}U_1(t)+\mu & n=1\\ \sum_{k=0}^{n}\frac{1}{k!}\mu^k U_{n-k}(t)-\sum_{k=0}^{n-2}\frac{1}{k!}\mu^k U_{n-2-k}(t) & n\geq 2.\end{cases} \tag{3.8}$$

**Proof**. For the cuts defined by (2.7), we consider the following test function

$$\Psi_{\mu,n}(u,z)=e^{\mu(u-z+\sqrt{z^2-1})}(z-\sqrt{z^2-1})^n,\ n\geq 0. \tag{3.9}$$

For $s\in(-1,1)$, we have

$$\Psi^{\pm}_{\mu,n}(u,s)=\lim_{\varepsilon>0,\varepsilon\to 0}\Psi_{\mu,n}(u,s\pm i\varepsilon)=e^{\mu(u-s\pm i\sqrt{1-s^2})}(s\mp i\sqrt{1-s^2})^n. \tag{3.10}$$

Following the similar derivation in the proof of Lemma 1 we obtain

$$
\begin{aligned}
0 &= \frac{1}{2\pi i}\int_{|z|=d,d>1}\frac{1}{t-z}\frac{\Psi_{\mu,n}(u,z)}{\sqrt{z^2-1}}dz \\
&= \lim_{\substack{\varepsilon>0\\ \varepsilon\to 0}}\frac{1}{2\pi i}\int_{-1}^{1}\frac{1}{t-(s+i\varepsilon)}\frac{\Psi_{\mu,n}(u,s+i\varepsilon)}{\sqrt{(s+i\varepsilon)^2-1}}ds-\lim_{\substack{\varepsilon>0\\ \varepsilon\to 0}}\frac{1}{2\pi i}\int_{-1}^{1}\frac{1}{t-(s-i\varepsilon)}\frac{\Psi_{\mu,n}(u,s-i\varepsilon)}{\sqrt{(s-i\varepsilon)^2-1}}ds \\
&= \frac{1}{2}\frac{\Psi^{+}_{\mu,n}(u,t)-\Psi^{-}_{\mu,n}(u,t)}{i\sqrt{1-t^2}}+\frac{1}{2\pi i}\int_{-1}^{1}\frac{1}{t-s}\frac{\Psi^{+}_{\mu,n}(u,s)+\Psi^{-}_{\mu,n}(u,s)}{i\sqrt{1-s^2}}ds.
\end{aligned}
\tag{3.11}
$$

Reorganizing the terms in (3.11), for $n\geq 1$ we obtain

$$
\frac{1}{\pi}\int_{-1}^{1}\frac{\exp[\mu(t-s)]}{t-s}\frac{T_{\mu,n}(s)}{\sqrt{1-s^2}}ds=-U_{\mu,n-1}(t)\,. \tag{3.12}
$$

Notice that $\Psi_{\mu,n}(u,\infty)=0$ for $n\geq 1$, we have

$$
\begin{aligned}
0 &= \frac{1}{2\pi i}\int_{|z|=d,d>1}\frac{1}{t-z}\Psi_{\mu,n}(u,z)dz \\
&= \lim_{\substack{\varepsilon>0\\ \varepsilon\to 0}}\frac{1}{2\pi i}\int_{-1}^{1}\frac{1}{t-(s+i\varepsilon)}\Psi_{\mu,n}(u,s+i\varepsilon)ds-\lim_{\substack{\varepsilon>0\\ \varepsilon\to 0}}\frac{1}{2\pi i}\int_{-1}^{1}\frac{1}{t-(s-i\varepsilon)}\Psi_{\mu,n}(u,s-i\varepsilon)ds\,. \\
&= \frac{1}{2}[\Psi^{+}_{\mu,n}(u,t)+\Psi^{-}_{\mu,n}(u,t)]+\frac{1}{2\pi i}\int_{-1}^{1}\frac{1}{t-s}[\Psi^{+}_{\mu,n}(u,s)-\Psi^{-}_{\mu,n}(u,s)]ds.
\end{aligned}
\tag{3.13}
$$

It follows that for $n\geq 1$

$$
\frac{1}{\pi}\int_{-1}^{1}\frac{\exp[\mu(t-s)]}{t-s}U_{\mu,n-1}(s)\sqrt{1-s^2}\,ds=T_{\mu,n}(t)\,. \tag{3.14}
$$

We call (3.12) and (3.14) as the positive half of the chFHT. To find the formulas for the negative half, we consider the reciprocal of (3.9)

$$
\Psi_{-\mu,-n}(u,z)=\frac{1}{2}e^{-\mu(u-z+\sqrt{z^2-1})}(z-\sqrt{z^2-1})^{-n}\,,\ n\geq 0\,. \tag{3.15}
$$

If $z\to\infty$, $(z-\sqrt{z^2-1})^{-n}$ tends to infinity for $n\geq 1$ so that the preceding computations are not applicable. Let $x=z-\sqrt{z^2-1}$, and then $z=\frac{1}{2}(x+\frac{1}{x})$ is the Joukowsky transform. It follows that

$$
\begin{aligned}
&\frac{1}{2\pi i}\int_{|z|=d,d>1}\frac{1}{t-z}\frac{\Psi_{-\mu,-n}(u,z)}{\sqrt{z^2-1}}dz \\
&= \frac{1}{2\pi i}\int_{|x|=\delta,\delta<1}\frac{2x}{2xt-x^2-1}\frac{2x}{1-x^2}[\frac{1}{2}e^{-\mu(u-x)}\frac{1}{x^n}]\frac{1}{2x^2}(x^2-1)dx \\
&= \frac{1}{2\pi i}\int_{|x|=\delta,\delta<1}\frac{e^{-\mu u}e^{\mu x}}{1-2xt+x^2}\frac{1}{x^n}dx \\
&= \frac{1}{2\pi i}e^{-\mu u}\int_{|x|=\delta,\delta<1}\{\sum_{m=0}^{\infty}[\sum_{k=0}^{m}\frac{1}{k!}\mu^k U_{m-k}(t)]x^m\}\frac{1}{x^n}dx.
\end{aligned}
\tag{3.16}
$$

Here $1/(1-2xt+x^2)=\sum_{k=0}U_k(t)x^k$ has been used. The residue theorem yields

$$
\frac{1}{2\pi i}\int_{|z|=d,d>1}\frac{1}{t-z}\frac{\Psi_{-\mu,-n}(u,z)}{\sqrt{z^2-1}}dz=e^{-\mu u}\sum_{k=0}^{n-1}\frac{1}{k!}\mu^k U_{n-1-k}(t)\,. \tag{3.17}
$$

Around $[-1,\ 1]$ we have the following identities

$$\begin{aligned}
&\frac{1}{2\pi i}\int_{|z|=d,d>1}\frac{1}{t-z}\frac{\Psi_{-\mu,-n}(u,z)}{\sqrt{z^2-1}}dz \\
&=\lim_{\substack{\varepsilon>0\\ \varepsilon\to 0}}\frac{1}{2\pi i}\int_{-1}^{1}\frac{1}{t-(s+i\varepsilon)}\frac{\Psi_{-\mu,-n}(u,s+i\varepsilon)}{\sqrt{(s+i\varepsilon)^2-1}}ds-\lim_{\substack{\varepsilon>0\\ \varepsilon\to 0}}\frac{1}{2\pi i}\int_{-1}^{1}\frac{1}{t-(s-i\varepsilon)}\frac{\Psi_{-\mu,-n}(u,s-i\varepsilon)}{\sqrt{(s-i\varepsilon)^2-1}}ds \\
&=\frac{1}{2}\frac{\Psi^{+}_{-\mu,-n}(u,t)-\Psi^{-}_{-\mu,-n}(u,t)}{i\sqrt{1-t^2}}+\frac{1}{2\pi i}\int_{-1}^{1}\frac{1}{t-s}\frac{\Psi^{+}_{-\mu,-n}(u,s)+\Psi^{-}_{-\mu,-n}(u,s)}{i\sqrt{1-s^2}}ds.
\end{aligned} \tag{3.18}$$

Reorganizing the terms in (3.18), for $n\ge 1$ we obtain

$$\frac{1}{\pi}\int_{-1}^{1}\frac{\exp[-\mu(t-s)]}{t-s}\frac{T_{\mu,n}(s)}{\sqrt{1-s^2}}ds=U_{\mu,n-1}(t)-2e^{-\mu u}\sum_{k=0}^{n-1}\frac{1}{k!}\mu^k U_{n-1-k}(t)\,. \tag{3.19}$$

Similarly, for $n\ge 1$ we have

$$\begin{aligned}
&\frac{1}{2\pi i}\int_{|z|=d,d>1}\frac{1}{t-z}\Psi_{-\mu,-n}(u,z)dz \\
&=\frac{1}{2\pi i}\int_{|x|=\delta,\delta<1}\frac{2x}{2xt-x^2-1}\Big[\frac{1}{2}e^{-\mu(u-x)}\frac{1}{x^n}\Big]\frac{1}{2x^2}(x^2-1)dx \\
&=\frac{1}{2\pi i}\int_{|x|=\delta,\delta<1}\frac{e^{-\mu u}e^{\mu x}}{1-2xt+x^2}\frac{1-x^2}{2x^{n+1}}dx \\
&=\frac{1}{2\pi i}\int_{|x|=\delta,\delta<1}\frac{e^{-\mu u}}{2}\Big\{\sum_{m=0}^{\infty}\Big[\sum_{k=0}^{m}\frac{1}{k!}\mu^k U_{m-k}(t)\Big]x^m\Big\}\frac{1-x^2}{x^{n+1}}dx.
\end{aligned} \tag{3.20}$$

By the residue theorem for $n\ge 1$ we have

$$\frac{1}{2\pi i}\int_{|z|=d,d>1}\frac{1}{t-z}\Psi_{-\mu,-n}(u,z)dz=\frac{e^{-\mu t}}{2}\begin{cases}U_1(t)+\mu & n=1\\ \sum_{k=0}^{n}\frac{1}{k!}\mu^k U_{n-k}(t)-\sum_{k=0}^{n-2}\frac{1}{k!}\mu^k U_{n-2-k}(t) & n\ge 2.\end{cases} \tag{3.21}$$

Then for $n\ge 1$, we obtain

$$\begin{aligned}
&\frac{1}{\pi}\int_{-1}^{1}\frac{\exp[-\mu(t-s)]}{t-s}U_{\mu,n-1}(s)\sqrt{1-s^2}ds \\
&\quad=-T_{\mu,n}(t)+e^{-\mu t}\begin{cases}U_1(t)+\mu & n=1\\ \sum_{k=0}^{n}\frac{1}{k!}\mu^k U_{n-k}(t)-\sum_{k=0}^{n-2}\frac{1}{k!}\mu^k U_{n-2-k}(t) & n\ge 2.\end{cases}
\end{aligned} \tag{3.22}$$

Identities (3.17) and (3.20) are the negative half of the chFHT. Combining the positive and negative halves we complete the proof.

**Remark 4**. It is interesting to notice that (3.7) is similar to (3.38) of [20]. If $\mu=0$, (3.7) and (3.8) degenerate to (3.1) and (3.2), respectively due to the identity $T_n(t)=0.5(U_n(t)-U_{n-2}(t))$. For the cases of $n$=1 and 2, we performed some numerical experiments to verify (3.8). Also, it is straightforward by (3.7) and (3.8) to derive the chFHTs of following functions

$$\frac{\cos(\mu\sqrt{1-s^2})T_n(s)}{\sqrt{1-s^2}}=\frac{T_{\mu,n}(s)+T_{-\mu,n}(s)}{2}\frac{1}{\sqrt{1-s^2}}\,, \tag{3.23}$$

$$\cos(\mu\sqrt{1-s^2})U_n(s)\sqrt{1-s^2}=\frac{U_{\mu,n}(s)+U_{-\mu,n}(s)}{2}\sqrt{1-s^2}\,. \tag{3.24}$$

# IV. Numerical experiments

In practical applications, the source function to be inversed is always bounded and continuous; thus $f(t)\in L_d^2(I)$ and (2.27) is preferred to avoid $\int_{-1}^{1}\cosh(\mu t)f(t)dt$. The discretization of (2.27) can be realized by the Chebyshev series for fast and stable numerical implementation. Here we describe the idea to convert singular integrals at two ends in (2.27) to normal integrals. Define the coordinate transformation

$$s=\cos\theta,\ t=\cos\varphi,\qquad \theta,\ \varphi\in(0,\pi). \tag{4.1}$$

Under (4.1), the inverse ${}^{d}H_{\mu}^{-1}F_{\mu}$ of (2.27) becomes

$$\begin{aligned}{}^{d}H_{\mu}^{-1}F_{\mu}(\cos\varphi)&=\cos(\mu\sin\varphi)\sin\varphi\frac{1}{\pi}\int_{0}^{\pi}\frac{\cos(\mu\sin\theta)}{\cos\theta-\cos\varphi}F_{\mu}(\cos\theta)d\theta\\&+\sin(\mu\sin\varphi)\frac{1}{\pi}\int_{0}^{\pi}\frac{\sin(\mu\sin\theta)\sin\theta}{\cos\theta-\cos\varphi}F_{\mu}(\cos\theta)d\theta.\end{aligned} \tag{4.2}$$

With (4.1), $L_d^2(I)$ becomes $L^2[0,\pi]$, and we use $\|{}^{d}H_{\mu}^{-1}F_{\mu}-f\|_{L^2[0,\pi]}$ to measure the deviation of ${}^{d}H_{\mu}^{-1}F_{\mu}$ from $f$. The value of $\|{}^{d}H_{\mu}^{-1}F_{\mu}-f\|_{L^2[0,\pi]}^{2}$ is commonly called the mean square error (MSE) in the data science. To measure the difference between the ground truth data and derived data, in this paper we use the *data-error-ratio* (DER) which is defined as

$$\text{DER}=\log_{10}\frac{\|\,data\,\|_{L^2[0,\pi]}}{\|\,error\,\|_{L^2[0,\pi]}}. \tag{4.3}$$

At the Chebyshev nodes, (4.2) can be implemented by the sine and cosine transforms. Sine and cosine transforms are orthogonal so that the numerical realizations are stable as shown from the early results in [18, 19] and the results in this paper. The Chebyshev nodes are

$$t_m=s_m=\cos(\frac{m+0.5}{N}\pi),\ m=0,\cdots,N-1. \tag{4.4}$$

Define two matrices $\mathbf{C}=\{c(m,n)\}$ and $\mathbf{S}=\{s(m,n)\}$, $m,n=0,\cdots,N-1$

$$c(m,n)=\sqrt{\frac{2}{N}}\cos(\frac{m+0.5}{N}n\pi),\ s(m,n)=\sqrt{\frac{2}{N}}\sin(\frac{m+0.5}{N}n\pi). \tag{4.5}$$

The Hilbert transform and its inverse can be realized by $\mathbf{H}^{-1}=\mathbf{S}\mathbf{C}^{T}$. Other discretization schemes were explored, but (4.4) and (4.5) produced the superior numerical stability for large constants such as $\mu=8\pi, 8\pi i$ and $20-20i$.

From (2.3), (2.5), (3.7) and (3.8), we know series of functions and their chFHTs. The goal of numerical experiments is to compare the source function $f$ and numerically inversed ${}^{d}H_{\mu}^{-1}F_{\mu}$. All the results were computed by the Python Numpy library. Total 1000 angular points were evenly sampled and displayed on Chebyshev nodes by (4.4). The display of DER is restricted to the precision with 2 decimal digits. In this section we present the numerical results for the pair defined by (2.5) from a variety of value ranges for $\mu$ with or without Gaussian white noise. Under the coordinate transformation of (4.1), the function pair of (2.5) is

$$\begin{cases}f(\cos\varphi)&=\cos(\mu\sin\varphi)\sin\varphi\\F_{\mu}(\cos\theta)&=\cos\theta\cosh(\mu\cos\theta)-\dfrac{\mu}{2}\sinh(\mu\cos\theta).\end{cases} \tag{4.6}$$

### A. Superior numerical stability for large constants

In the first numerical experiment, we took two real and two imaginary $\mu$ values of $4\pi$, $8\pi$, $4\pi i$ and $8\pi i$. All the plots are shown in Fig 1.

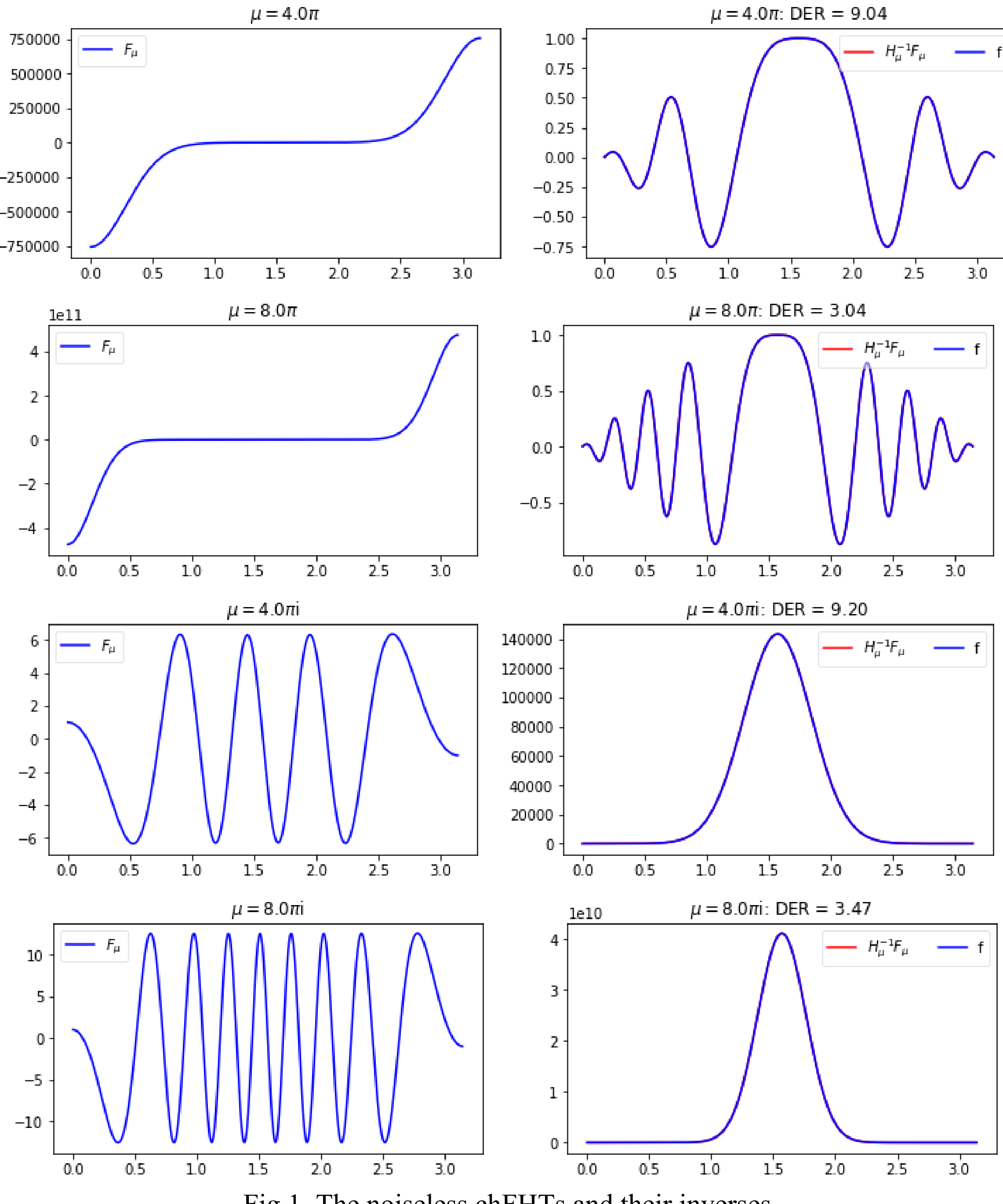

Fig 1. The noiseless chFHTs and their inverses

In Fig.1, the left side is $F_{\mu}$ and the right side is the overlay display of the source function $f$ and numerical inverse ${}^{d}H_{\mu}^{-1}F_{\mu}$. One observation from Fig. 1 is that the accuracy of the inverse is extremely high for large $\mu$ even for the value range up to $10^{10}$. This is something we were not be able to achieve with the Fourier analysis in [9] and the FBP algorithm in [17].

### B. Moderate constant with noise

In the second numerical experiment, we added Gaussian white noise $n(s,\sigma)$ to $F_\mu(s)$, here $n(s,\sigma)$ is the Gaussian random numbers with zero mean and standard deviation of $\sigma$. Let $\widetilde{F}_\mu(s) = F_\mu(s) + n(s,\sigma)$ and denote by $\widetilde{f}(t) = {}^{d}H_\mu^{-1}\widetilde{F}_\mu(t)$ the inverse from the noisy $\widetilde{F}_\mu(s)$. As expected, in the presence of noise, the inverse becomes worse for large $\mu$. In this experiment, we take $\mu = \pi$, $2\pi$, $\pi i$ and $2\pi i$, respectively, present the results of the inverse with DERs for $\sigma = 0.1$ in Fig. 2.

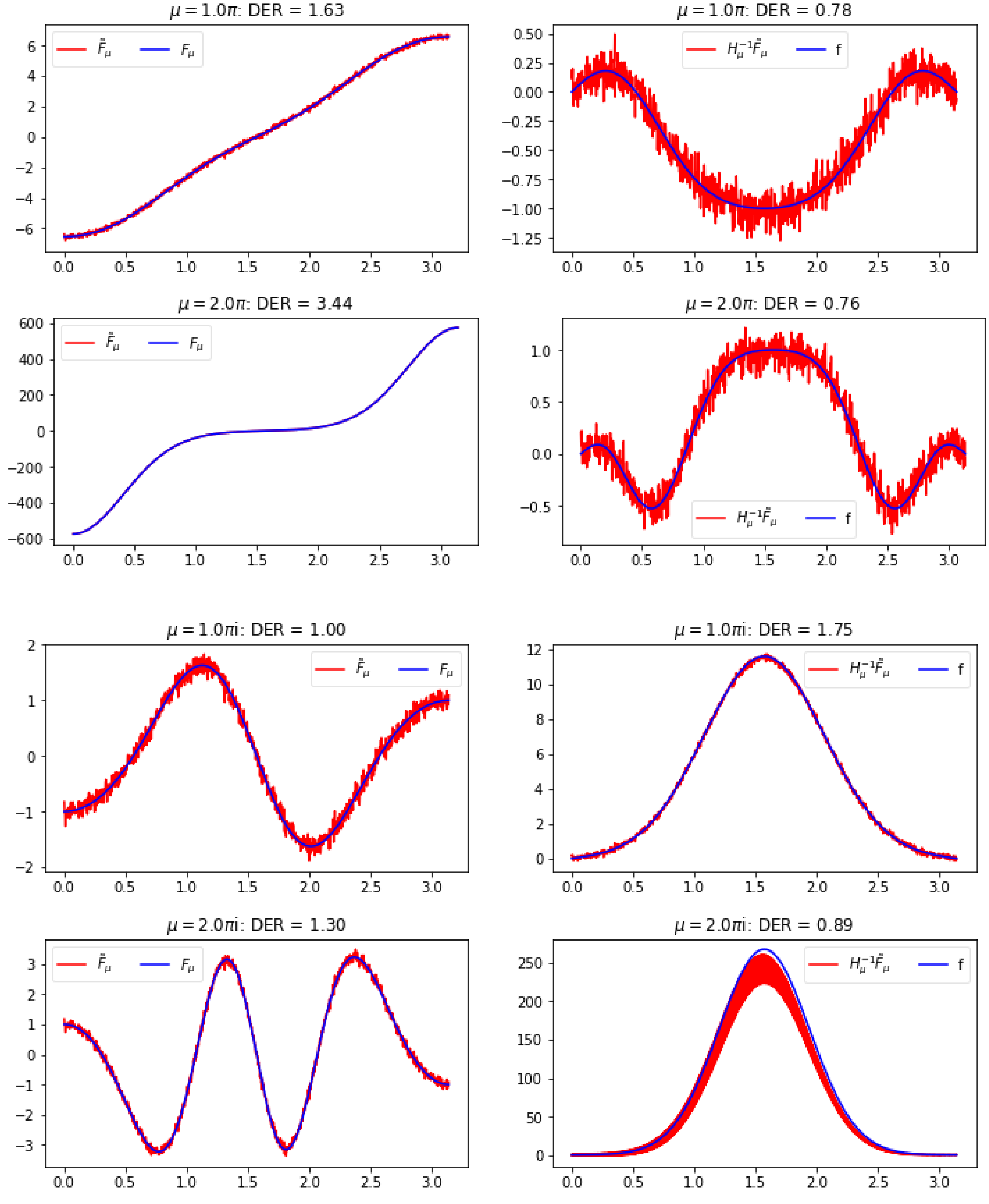

Fig. 2 The noisy chFHTs and their inverses

### C. Superior numerical stability for large complex constants

In the third numerical experiment for complex $\mu$, we took $\mu = 10.0 + 10.0i$ and $\mu = 20.0 - 20.0i$. The results are shown in Fig. 3.

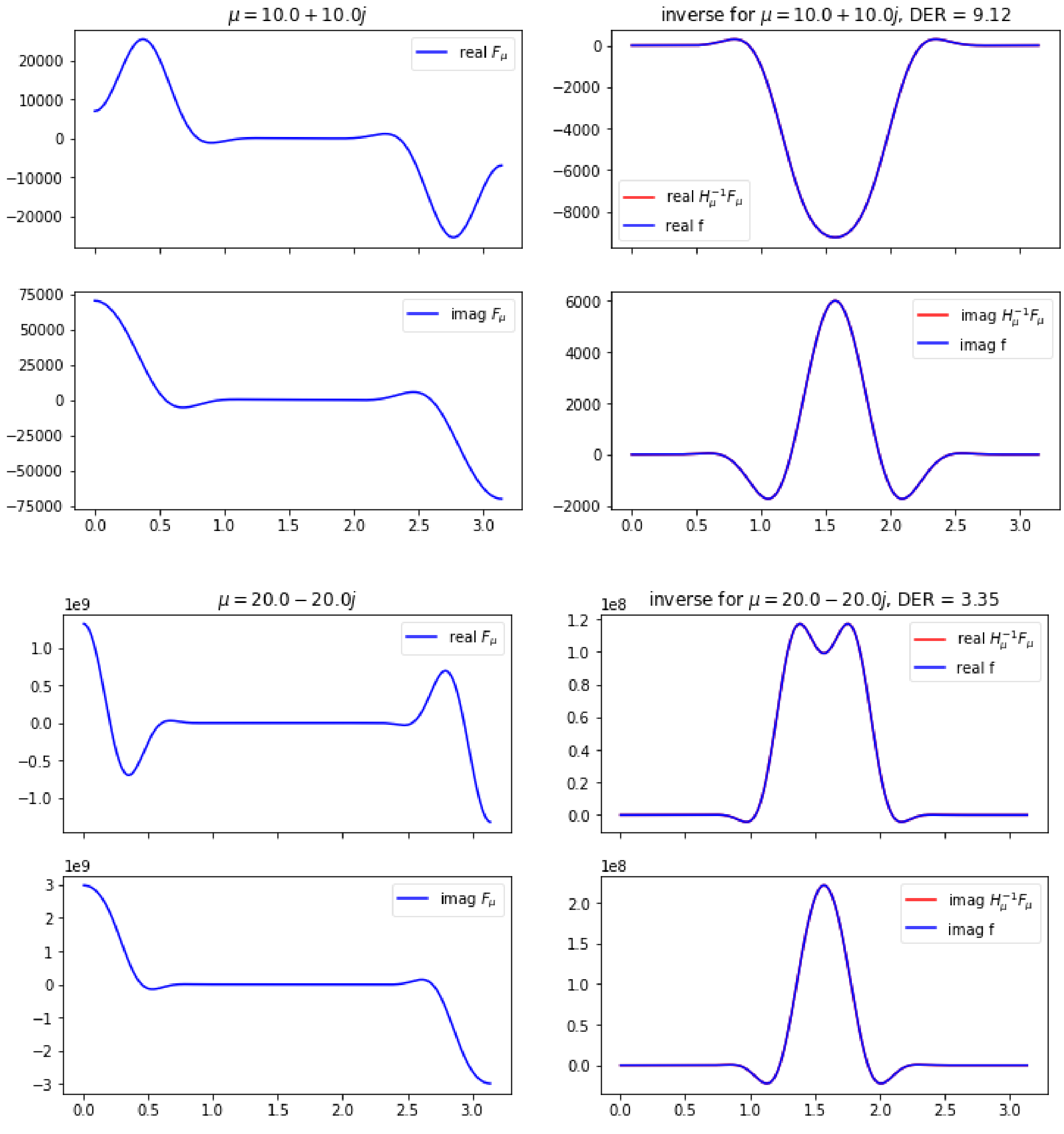

Fig. 3 noiseless chFHTs and their inverse from complex constants

For both $\mu = 10.0 + 10.0i$ and $\mu = 20.0 + 20.0i$, the accuracy of the inverse is unusually high for the value range up to $10^9$. In our study, sine and cosine transforms of (4.5) at the Chebyshev nodes are the main factor leading to such high accuracy. For example, a direct discretization of (2.27) at the evenly-sampled points in $[-1, 1]$ is not able to produce such level of accuracy.

### D. Moderate complex constants with noise

In the fourth numerical experiment, we added Gaussian white noise for constants $\mu = 2.0 + 2.0i$ and $4.0 - 4.0i$. The results are shown in Fig. 4.

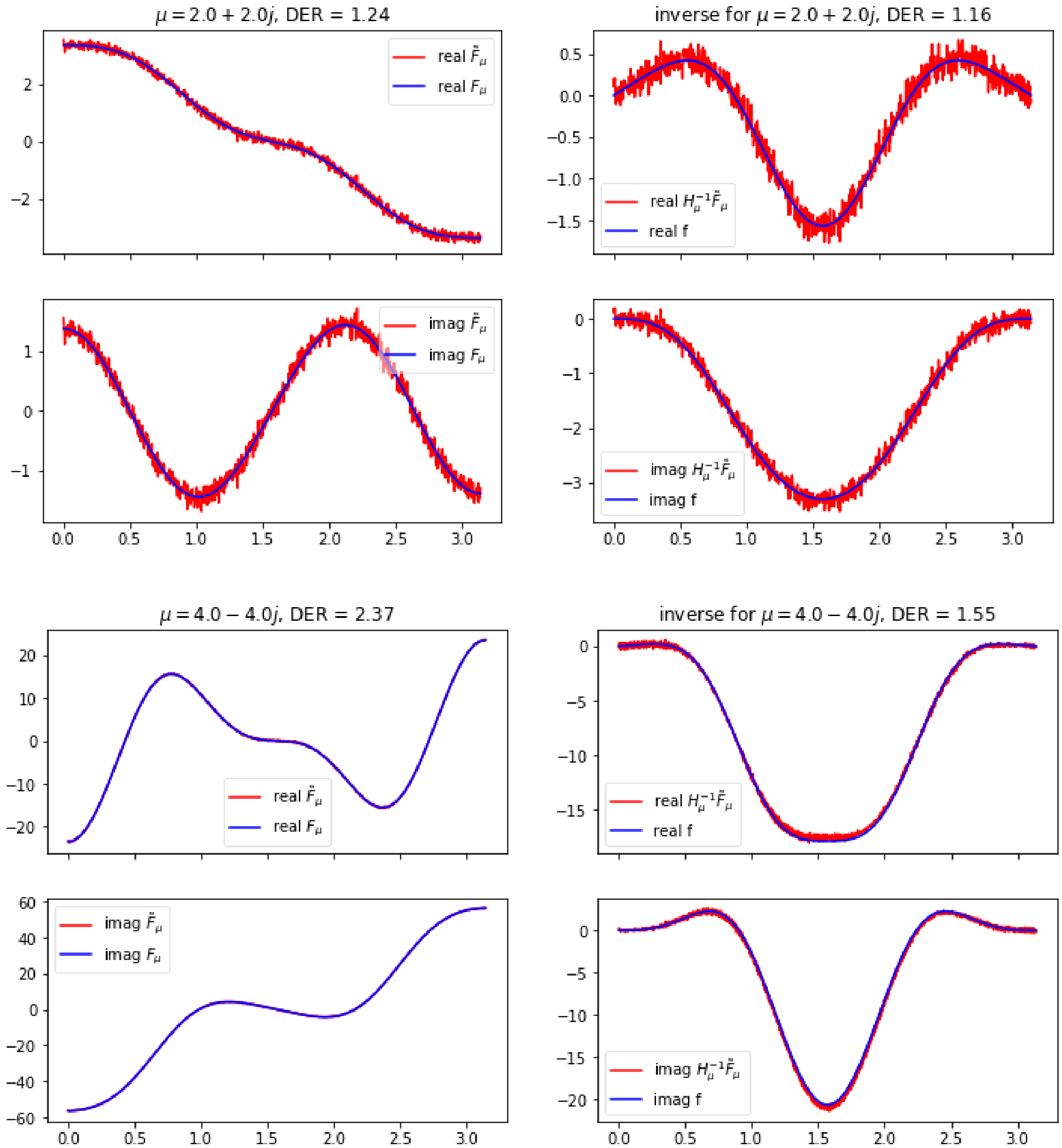

Fig. 4 The noisy chFHTs and their inverses from complex constants

The Python implementation of (2.27) was straightforward and has been posted on GitHub for readers to study the numerical characteristics of (2.27). For example, it is interesting to observe the different numerical characteristics from other pairs defined by (3.8)

$$\begin{cases} f(\cos\varphi) & = \sin(2\varphi - \mu\sin\varphi) \\ F_\mu(\cos\theta) & = 0.5\exp(\mu\cos\theta)[4\cos^2\theta - 2 + 2\mu\cos\theta + 0.5\mu^2]. \end{cases} \tag{4.7}$$

Here (4.7) is (3.8) for $n = 2$. Our Python implementation includes more test functions. All the numerical experiments don't have any data characteristics pertinent to particular applications. It is interesting to study the inversion characteristics in SPECT and other type of tomographic imaging with application dependent noise model other than Gaussian white noise.

## Discussion and conclusions

The study on the inversion of the chFHT was initiated by [14]. With the Bertola-Katsevich-Tovbis inversion formulas [2], it seems that the theoretical study reaches a point at which the half scan formulated in [14] is solved with the explicit inversion formulas (2.27) and (2.28). The constant $\mu$ can be complex in (2.27) and (2.28), thus (1.2) is also useful to ultrasonic tomography and Doppler tomography imaging described in [9, 13, 17]. The explicit expressions of the chFHT of the exponential Chebyshev functions may be useful to construct the basis in defining a cost function to consider the inversion from partial data. It is known that $T_n(t)$ and $U_n(t)$ construct a complete orthogonal basis in $L_d^2(I)$ and $L_m^2(I)$, respectively. However, if $\mu \neq 0$, $\{T_{\mu,n}(t)\}$ and $\{U_{\mu,n}(t)\}$ are not orthogonal and it is unknown if they can construct a complete basis in $L_d^2(I)$ and $L_m^2(I)$.

Initial numerical results showed the feasibility of (2.27) for practical use, certainly more quantitative numerical analysis from complicated data pertinent to particular applications is desirable. Here it is worth mentioning the advantage of (2.27) over (2.28). The weighted integral $\int_{-1}^{1}\cosh(\mu t) f(t)dt$ is required in (2.28) so that any noise or perturbation to that weighted integral could lead to deviations to entire function $f(t)$. From application point of view, the confidence of using (2.28) is unpredictable in the presence of noise and error. The factor $1/\sqrt{1-t^2}$ could significantly increase the errors around the two endpoints so that (2.28) may not be robust to errors either. Fortunately, formula (4.2) as the angular version of (2.27) can void the instability of the weighted integral and the error amplification at two endpoints. Another important numerical realization technique is the cosine and sine transforms at the Chebyshev nodes. The direct discretization of the Hilbert transform on the unit interval does not generate as high accuracy in our numerical experiments as the angular version (4.2). In the presence of noise, the error increases significantly for a large $\mu$, this may require noise treatment in applications with a large $\mu$.

*Python implementation by Jupyter Notebook can be found from GitHub*:

*https://github.com/jshyou/Medical-Imaging/blob/master/complexCoshFHT.html*